\documentclass{gtart_h}  

\def\ifplaintex{\expandafter\ifx\csname documentclass\endcsname\relax}

\ifplaintex 
\else
\fi

\expandafter\ifx\csname beginpicture\endcsname\relax
\expandafter\ifx\csname documentclass\endcsname\relax
\input pictex \else
\input prepictex \input pictex \input postpictex \fi\fi

\def\gt{{\mathsurround=0pt\it $\cal G\mskip-2mu$eometry \&\ 
$\cal T\!\!$opology}}        

\def\gtp{{\mathsurround=0pt\it $\cal G\mskip-2mu$eometry \&\ 
$\cal T\!\!$opology $\cal P\!$ublications}}  

\def\lognumber#1{\def\thelognumber{#1}}
\def\volumenumber#1{\def\thevolumenumber{#1}}
\def\papernumber#1{\def\thepapernumber{#1}}
\def\volumeyear#1{\def\thevolumeyear{#1}}

\def\pagenumbers#1#2{\def\startpage{#1}\def\finishpage{#2}}
\def\published#1{\def\publishdate{#1}}
\def\proposed#1{\def\theproposer{#1}}
\def\seconded#1{\def\theseconders{#1}}
\def\received#1{\def\receiveddate{#1}}
\def\revised#1{\def\reviseddate{#1}}
\def\accepted#1{\def\accepteddate{#1}}

\def\asciiaddress#1{\def\theasciiaddress{#1}}

\long\def\asciiabstract#1{\long\def\theasciiabstract{#1}}

\let\\\par\let\thelognumber\relax
\let\thevolumenumber\relax\let\thepapernumber\relax
\let\thevolumeyear\relax\let\thesamplenumber\relax\let\startpage\relax
\let\finishpage\relax\let\publishdate\relax\let\receiveddate\relax
\let\reviseddate\relax\let\accepteddate\relax\let\theasciititle\relax
\let\theasciiauthors\relax\let\theasciiaddress\relax
\let\theasciiabstract\relax
\let\theasciiemail\relax\let\theshortauthors\relax\let\theshorttitle\relax

\long\def\maketitlep{   

\count0=\startpage

\gt\hfill      
\beginpicture
\setcoordinatesystem units <0.33truein, 0.33truein> point at 2.2 0.9
\setplotsymbol ({$\cal G$})
\plotsymbolspacing=9truept
\circulararc 315 degrees from 0 1 center at 0 0
\setplotsymbol ({$\cal T$})
\circulararc 315 degrees from 1 -1 center at 1 0
\endpicture
\break
{\small\ifx\thesamplenumber\relax 
Volume \else Sample
\fi\thevolumenumber\ (\thevolumeyear)
\startpage--\finishpage\nl
Published: \publishdate}
\vglue 0.5truein plus 0.4fil minus 0.1truein

{\parskip=0pt\leftskip 0pt plus 1fil\def\\{\par\smallskip}{\ifplaintex\large
\else\Large\fi\bf\thetitle}\par\medskip}   

\vglue 0pt plus 0.1fil 

{\parskip=0pt\leftskip 0pt plus 1fil\def\\{\par}{\sc\theauthors}
\par\medskip}

\vglue 0pt plus 0.1fil 

{\small\parskip=0pt\let\newline\\
{\leftskip 0pt plus 1fil\def\\{\par}{\sl\theaddress}\par}
\expandafter\ifx\theemail\relax    
\relax\else\vglue 5pt plus 0.02fil minus 2pt\def\\{\stdspace{\rm 
and}\stdspace} 
\cl{Email:\stdspace\tt\theemail}\fi
\ifx\theurl\relax                  
\relax\else\vglue 5pt plus 0.02fil minus 2pt\def\\{\stdspace{\rm 
and}\stdspace}
\cl{URL:\stdspace\tt\theurl}\fi\par}

\vglue 7pt plus 0.3fil minus 3pt

{\bf Abstract}
\vglue 5pt plus 0.1fil minus 2pt

\theabstract

\vglue 7pt plus 0.3fil minus 3pt

{\bf AMS Classification numbers}\quad Primary:\quad \theprimaryclass

Secondary:\quad \thesecondaryclass

\vglue 5pt plus 0.3fil minus 2pt

{\bf Keywords:}\quad \thekeywords

\vglue 10pt plus 0.5fil minus 5pt

{\small  Proposed: \theproposer\hfill Received: \receiveddate\nl
Seconded: \theseconders\hfill 
\ifx\reviseddate\relax                         
Accepted: \accepteddate                        
\else
Revised: \reviseddate                          
\fi}
\eject
}       

\let\maketitlepage\maketitlep
\let\maketitle\maketitlepage

\font\phead=cmsl9 scaled 950
\font=cmsl9 scaled 1050
\font\pnum=cmbx10 scaled 913
\font\lnum=cmbx10 
\font\pfoot=cmsl9 scaled 950
\font=cmsl9 scaled 1050
\ifplaintex
\headline{\vbox to 0pt{\vskip -4.5mm\line{\small\phead\ifnum
\count0=\startpage ISSN 1364-0380 (on line)
1465-3060 (printed) \hfill {\pnum\folio}\else\ifodd\count0\def\\{ }%
\ifx\theshorttitle\relax\thetitle\else\theshorttitle\fi\hfill{\pnum\folio}
\else\def\\{ and }{\pnum\folio}\hfill\ifx\theshortauthors\relax\theauthors
\else\theshortauthors\fi\fi\fi}\vss}}
\footline{\vbox to 0pt{\vglue 0mm\line{\small\pfoot\ifnum\count0=\startpage
\copyright\ \gtp\hfill\else
\gt, Volume \thevolumenumber\ (\thevolumeyear)\hfill\fi}\vss
}}
\else
\makeatletter
\def\@oddhead{{\small\lhead\ifnum\count0=\startpage ISSN 1364-0380 (on line)
1465-3060 (printed) \hfill {\lnum\number\count0}\else\ifodd\count0
\def\\{ }\ifx\theshorttitle\relax \thetitle \else\theshorttitle\fi\hfill
{\lnum\number\count0}\else\def\\{ and }{\lnum\number\count0}
\hfill\ifx\theshortauthors\relax 
\theauthors\else\theshortauthors\fi\fi\fi}}\def\@evenhead{\@oddhead}
\def\@oddfoot{\small\lfoot\ifnum\count0=\startpage\copyright\ \gtp\hfill\else
\gt, Volume \thevolumenumber\ (\thevolumeyear)\hfill\fi}
\def\@evenfoot{\@oddfoot}
\makeatother
\fi

\newwrite\gtoutfile
\long\gdef\makeheadfile{  
{\def\\{, }\def\s{ }
\immediate\openout\gtoutfile head.xxx
\immediate\write\gtoutfile{Proxy-for: \ifx\theasciiauthors\relax
\theauthors\else\theasciiauthors\fi\s<\ifx\theasciiemail\relax\theemail\else\theasciiemail\fi>}
\immediate\write\gtoutfile{\noexpand\\}
\immediate\write\gtoutfile{Authors: \ifx\theasciiauthors\relax
\theauthors\else\theasciiauthors\fi}
{\def\\{ }\immediate\write\gtoutfile{Title: \ifx\theasciititle\relax
\thetitle\else\theasciititle\fi}}
\immediate\write\gtoutfile{Subj-class: GT or SG or MG etc}
\immediate\write\gtoutfile{MSC-class: \theprimaryclass\ifx\thesecondaryclass\relax\else, \thesecondaryclass\fi}
\immediate\write\gtoutfile{Journal-ref: Geom. Topol. \thevolumenumber
(\thevolumeyear) \startpage-\finishpage}
\immediate\write\gtoutfile{Comments: Published by Geometry and Topology at}
\immediate\write\gtoutfile{\s\s http://www.maths.warwick.ac.uk/gt/GTVol\thevolumenumber/paper\thepapernumber.abs.html}
\immediate\write\gtoutfile{\noexpand\\}
\immediate\write\gtoutfile{}
\ifx\theasciiabstract\relax
\immediate\write\gtoutfile{\theabstract}\else
\immediate\write\gtoutfile{\theasciiabstract}\fi
\immediate\write\gtoutfile{}
\immediate\write\gtoutfile{\noexpand\\}
\immediate\write\gtoutfile{}
\immediate\closeout\gtoutfile}}  

\def\maketitlepage{\maketitlep\makeheadfile}
\let\maketitle\maketitlepage

\lognumber{404}
\received{6 January 2004}
\volumenumber{8}\papernumber{35}\volumeyear{2004}
\pagenumbers{1281}{1300}
\revised{21 September 2004}
\published{28 September 2004}
\accepted{21 September 2004}
\proposed{Joan Birman}
\seconded{Robion Kirby, Cameron Gordon}

\usepackage{amsmath,amssymb}
\let\Bbb\mathbb

\newtheorem{thm}{Theorem}[section]
\newtheorem{lem}[thm]{Lemma}
\newtheorem{prop}[thm]{Proposition}
\newtheorem{cor}[thm]{Corollary}

\theoremstyle{definition}

\newtheorem*{acknow}{Acknowledgement}
\newtheorem*{rmk}{Remarks}
\newtheorem*{case1}{Case 1}
\newtheorem*{case2}{Case 2}
\newtheorem*{case3}{Case 3}

\begin{document}

\def\C{{\mathbb C}}\def\Z{{\mathbb Z}}\def\MM{{\mathcal M}}
\def\Id{{\rm Id}}\def\deg{{\rm deg}}\def\N{{\mathbb N}}
\def\Sym{{\rm Sym}}\def\AA{{\mathcal A}}\def\UU{{\mathcal U}}
\def\lg{{\rm lg}}

\title{The proof of Birman's conjecture on\\singular braid monoids}
\authors{Luis Paris}
\address{Institut de Math\'ematiques de Bourgogne, Universit\'e de 
Bourgogne\\UMR 5584 du CNRS, BP 47870, 21078 Dijon cedex,
France}
\asciiaddress{Institut de Mathematiques de Bourgogne, Universite de 
Bourgogne\\UMR 5584 du CNRS, BP 47870, 21078 Dijon cedex,
France}
\email{lparis@u-bourgogne.fr}
\urladdr{http://math.u-bourgogne.fr/topo/paris/index.html}

\begin{abstract}
Let $B_n$ be the Artin braid group on $n$ strings with standard
generators $\sigma_1, \dots,$ $\sigma_{n-1}$, and let $SB_n$ be the
singular braid monoid with generators $\sigma_1^{\pm 1}, \dots,
\sigma_{n-1}^{\pm 1},$ $\tau_1, \dots, \tau_{n-1}$. The
desingularization map is the multiplicative homomorphism $\eta\co SB_n
\to \Bbb Z[B_n]$ defined by $\eta(\sigma_i^{\pm 1}) = \sigma_i^{\pm
1}$ and $\eta (\tau_i) = \sigma_i - \sigma_i^{-1}$, for $1 \le i \le
n-1$.  The purpose of the present paper is to prove Birman's
conjecture, namely, that the desingularization map $\eta$ is
injective.
\end{abstract}
\asciiabstract{%
Let B_n be the Artin braid group on n strings with standard generators
sigma_1, ..., sigma_{n-1}, and let SB_n be the singular braid monoid
with generators sigma_1^{+-1}, ..., sigma_{n-1}^{+-1}, tau_1, ...,
tau_{n-1}. The desingularization map is the multiplicative
homomorphism eta: SB_n --> Z[B_n] defined by eta(sigma_i^{+-1}) =
_i^{+-1} and eta(tau_i) = sigma_i - sigma_i^{-1}, for 1 <= i <= n-1.
The purpose of the present paper is to prove Birman's conjecture,
namely, that the desingularization map eta is injective.}


\maketitlepage

\section{Introduction}

Define an {\it $n$--braid} to be a collection $b=(b_1, \dots, b_n)$ of disjoint smooth paths in 
$\C \times [0,1]$, called the {\it strings} of $b$, such that the $k$-th string $b_k$ runs 
monotonically in $t \in [0,1]$ from the point $(k,0)$ to some point $(\zeta(k),1)$, 
where $\zeta$ is a permutation of $\{1,2, \dots, n\}$.
An {\it isotopy} in this context is a deformation through braids which fixes the 
ends. Multiplication of braids is defined by concatenation. The isotopy classes of braids with 
this multiplication form a group, called {\it braid group on $n$ strings}, and denoted by $B_n$. 
This group has a well-known presentation with generators $\sigma_1, \dots, \sigma_{n-1}$ and 
relations
\[
\begin{array}{cl}
\sigma_j \sigma_k=\sigma_k\sigma_j\quad&\text{if}\qua |j-k|>1\,,\\
\sigma_j\sigma_k\sigma_j=\sigma_k\sigma_j\sigma_k\quad&\text{if}\qua |j-k|=1\,.
\end{array}
\]
The group $B_n$ has other equivalent descriptions as a group of automorphisms of a free group, as 
the fundamental group of a configuration space, or as the mapping class group of the $n$--punctured 
disk, and plays a prominent r\^ole in many disciplines. We refer to \cite{Bir1} for a general 
exposition on the subject.

The Artin braid group $B_n$ has been extended to the {\it singular braid monoid} $SB_n$ by Birman 
\cite{Bir2} and Baez \cite{Bae} in order to study Vassiliev invariants. The strings of a singular braid are 
allowed to intersect transversely in finitely many double points, called {\it singular points}. 
As with braids, isotopy is a deformation through singular braids which fixes the ends, and 
multiplication is by concatenation. Note that the isotopy classes of singular braids form a 
monoid and not a group. It is shown in \cite{Bir2} that $SB_n$ has a monoid presentation with 
generators $\sigma_1^{\pm 1}, \dots, \sigma_{n-1}^{\pm 1}$, $\tau_1, \dots, \tau_{n-1}$, and 
relations
\[
\begin{array}{cl}
\sigma_i \sigma_i^{-1} = \sigma_i^{-1} \sigma_i = 1\,, 
\quad \sigma_i \tau_i = \tau_i \sigma_i \quad
&\text{if}\qua 1 \le i\le n-1\,,\\
\sigma_i \sigma_j = \sigma_j \sigma_i\,, \quad 
\sigma_i \tau_j = \tau_j \sigma_i\,, \quad 
\tau_i \tau_j = \tau_j \tau_i\,,\quad
&\text{if}\qua |i-j|>1\,,\\
\sigma_i \sigma_j \sigma_i = \sigma_j \sigma_i \sigma_j\,, \quad 
\sigma_i \sigma_j \tau_i = \tau_j \sigma_i \sigma_j\,,\quad 
&\text{if}\qua |i-j|=1\,.
\end{array}
\] 
Consider the braid group ring $\Z[B_n]$. The natural embedding $B_n \to \Z[B_n]$ can be 
extended to a multiplicative homomorphism $\eta\co  SB_n \to \Z[B_n]$, called {\it 
desingularization map}, and defined by
\[
\eta (\sigma_i^{\pm 1}) = \sigma_i^{\pm 1}\,, \quad \eta(\tau_i) = \sigma_i-\sigma_i^{-1}\,, 
\quad \text{if}\qua 1\le i\le n-1\,.
\]
This homomorphism is one of the main ingredients of the definition of Vassiliev invariants for 
braids. It has been also used by Birman \cite{Bir2} to establish a relation between Vassiliev knot 
invariants and quantum groups.

One of the most popular problems in the subject, known as ``Birman's conjecture'', is to 
determine whether $\eta$ is an embedding (see \cite{Bir2}). 
At the time of this writing, the only 
known partial answer to this question is that $\eta$ is injective on singular braids with up to 
three singularities (see \cite{Zhu}), and on singular braids with up to three strings (see \cite{Jar}).

The aim of the present paper is to solve this problem, namely, we prove 
the following.

\begin{thm}
The desingularization map $\eta\co  SB_n \to \Z[B_n]$ is injective.
\end{thm}

Let $S_dB_n$ denote the set of isotopy classes of singular braids with $d$ singular points. 
Recall that a {\it Vassiliev invariant of type $d$} is defined to be a homomorphism $v\co  
\Z[B_n] \to A$ of $\Z$--modules which vanishes on $\eta( S_{d+1}B_n)$. 
One of the main results on Vassiliev braid invariants is that they
separate braids (see \cite{Bar}, \cite{Koh}, \cite{Pap}). 
Whether Vassiliev knot invariants separate knots remains an important open question.
Now, it has been shown by Zhu 
\cite{Zhu} that this separating property extends to singular braids if $\eta$ is injective. So, a 
consequence of Theorem 1.1 is the following.

\begin{cor}
Vassiliev braid invariants classify singular braids.
\end{cor}

Let $\Gamma$ be a graph (with no loop and no multiple edge), let $X$ be the set of vertices, 
and let $E=E(\Gamma)$ be the set of edges of $\Gamma$. Define the {\it graph monoid} of 
$\Gamma$ to be the monoid $\MM(\Gamma)$ given by the monoid presentation
\[
\MM(\Gamma)= \langle X\ |\ xy=yx\qua \text{if}\qua \{x,y\} \in E(\Gamma)\rangle^+\,.
\]
Graph monoids are also known as {\it free partially commutative monoids} or as {\it right-angled 
Artin monoids}. They were first introduced by Cartier and Foata \cite{CF} to study combinatorial 
problems on rearrangements of words, and, since then, have been extensively studied
by both computer scientists and mathematicians.

The key point of the proof of Theorem 1.1 consists in understanding the structure of the 
multiplicative submonoid of $\Bbb Z[B_n]$ generated by the set $\{ \alpha \sigma_i^2 \alpha^{-1} 
-1;\,\, \alpha \in B_n\ \text{and}\ 1\le i\le n-1\}$. More precisely, we prove the following.

\begin{thm}
Let $\Omega$ be the graph defined as follows.
\begin{itemize}
\item $\Upsilon=\{ \alpha \sigma_i^2 \alpha^{-1};\,\, \alpha \in B_n\ \text{and}\ 1\le i\le n-
1\}$ is the set of vertices of $\Omega$;
\item $\{u,v\}$ is an edge of $\Omega$ if and only if we have $uv=vu$ in $B_n$.
\end{itemize}
Let $\nu\co  \MM(\Omega) \to \Bbb Z[B_n]$ be the homomorphism defined by $\nu (u)=u-1$, for all 
$u\in\Upsilon$. Then $\nu$ is injective.
\end{thm}

The proof of the implication Theorem 1.3 $\Rightarrow$ Theorem 1.1 is based on the observation 
that $SB_n$ is isomorphic to the semi-direct product of $\MM(\Omega)$ with the braid group 
$B_n$, and that $\nu\co  \MM (\Omega) \to \Bbb Z [B_n]$ is the restriction to $\MM(\Omega)$ of 
the desingularization map. The proof of this implication is the subject of Section 2. Let 
$A_{i\,j}$, $1 \le i<j\le n$, be the standard generators of the pure braid group $PB_n$. In 
Section 3, we show that $\Upsilon$ is the disjoint union of the conjugacy classes of the 
$A_{i\,j}$'s in $PB_n$. Using homological arguments, we then show that we can restrict the 
study to the submonoid of $\MM(\Omega)$ generated by the conjugacy classes of two given 
generators, $A_{i\,j}$ and $A_{r\,s}$. If $\{i,j\} \cap \{r,s\} \neq \emptyset$, then the 
subgroup generated by the conjugacy classes of $A_{i\,j}$ and $A_{r\,s}$ is a free group, and we 
prove the injectivity using a sort of Magnus expansion (see Section 4). The case $\{i,j\} \cap 
\{r,s\} = \emptyset$ is handled using the previous case together with a technical result on automorphisms 
of free groups (Proposition 5.1).

\begin{acknow}
My first proof of Proposition 5.1 was awful, hence I asked some
experts whether they know another proof or a reference for the
result. The proof given here is a variant of a proof indicated to me
by Warren Dicks. So, I would like to thank him for his help.
\end{acknow}

\section{Theorem 1.3 implies Theorem 1.1}

We assume throughout this section that the result of Theorem 1.3 holds, and we prove Theorem 1.1.

Let $\delta_i= \sigma_i \tau_i$ for $1 \le i\le n-1$. Then $SB_n$ is generated as a monoid by 
$\sigma_1^{\pm 1}, \dots, \sigma_{n-1}^{\pm 1}$, $\delta_1, \dots, \delta_{n-1}$, and has a 
monoid presentation with relations
\[
\begin{array}{cl}
\sigma_i \sigma_i^{-1} = \sigma_i^{-1} \sigma_i = 1\,, 
\quad \sigma_i \delta_i = \delta_i \sigma_i\,, \quad 
&\text{if}\qua 1 \le i\le n-1\,,\\
\sigma_i \sigma_j = \sigma_j \sigma_i\,, \quad 
\sigma_i \delta_j = \delta_j \sigma_i\,, \quad 
\delta_i \delta_j = \delta_j \delta_i\,, \quad
&\text{if}\qua |i-j|>1\,,\\
\sigma_i \sigma_j \sigma_i = \sigma_j \sigma_i \sigma_j\,, \quad 
\sigma_i \sigma_j \delta_i = \delta_j \sigma_i \sigma_j\,, \quad 
&\text{if}\qua |i-j|=1\,.
\end{array}
\] 
Moreover, the desingularization map $\eta\co  SB_n \to \Z[B_n]$ is determined by
\[
\eta(\sigma_i^{\pm 1})= \sigma_i^{\pm 1}\,, \quad \eta(\delta_i)= \sigma_i^2-1\,, \quad 
\text{if}\qua 1\le i\le n-1\,.
\]
The following lemma is a particular case of \cite{FRZ}, Theorem 7.1.

\begin{lem}
Let $i,j \in \{ 1, \dots, n-1\}$, and let $\beta \in SB_n$. Then the 
following are equivalent:
\begin{enumerate}
\item $\beta \sigma_i^2 = \sigma_j^2 \beta$;
\item $\beta \delta_i = \delta_j \beta$.
\end{enumerate}
\end{lem}

This lemma shows the following.

\begin{lem}
Let $\hat \Omega$ be the graph defined as follows.
\begin{itemize}
\item $\hat \Upsilon= \{\alpha \delta_i \alpha^{-1};\,\, \alpha \in B_n\ \text{and}\ 1 \le i \le 
n-1\}$ is the set of vertices of $\hat \Omega$;
\item $\{\hat u, \hat v\}$ is an edge of $\hat \Omega$ if and only if we have $\hat u \hat v 
= \hat v \hat u$ in $SB_n$.
\end{itemize}
Then there exists an isomorphism $\varphi\co  \MM(\hat \Omega) \to \MM(\Omega)$ which sends 
$\alpha \delta_i \alpha^{-1}\in \hat \Upsilon$ to $\alpha \sigma_i^2 \alpha^{-1} \in \Upsilon$ 
for all $\alpha \in B_n$ and $1 \le i \le n-1$.
\end{lem}

\begin{proof}
Let $\alpha, \beta \in B_n$ and $i,j \in \{1, \dots, n-1\}$. Then, by Lemma 2.1,
\[
\begin{array}{cccc}
&\alpha \sigma_i^2 \alpha^{-1} = \beta \sigma_j^2 \beta^{-1}\quad
&\Leftrightarrow\quad
&(\beta^{-1} \alpha) \sigma_i^2 = \sigma_j^2 (\beta^{-1} \alpha)\\
\Leftrightarrow\quad
&(\beta^{-1} \alpha) \delta_i = \delta_j (\beta^{-1} \alpha)\quad
&\Leftrightarrow\quad
&\alpha \delta_i \alpha^{-1} = \beta \delta_j \beta^{-1}\,.
\end{array}
\]
This shows that there exists a bijection $\varphi\co  \hat \Upsilon \to \Upsilon$ which sends $\alpha 
\delta_i \alpha^{-1} \in \hat \Upsilon$ to $\alpha \sigma_i^2 \alpha^{-1}\in \Upsilon$ for all 
$\alpha \in B_n$ and $1 \le i\le n-1$. Let $\alpha, \beta \in B_n$ and $i,j \in \{1, \dots, n-
1\}$. Again, by Lemma 2.1,
\[
\begin{array}{cc}
&(\alpha \sigma_i^2 \alpha^{-1}) (\beta \sigma_j^2 \beta^{-1}) =
(\beta \sigma_j^2 \beta^{-1}) (\alpha \sigma_i^2 \alpha^{-1})\\
\Leftrightarrow\quad&
\sigma_i^2 (\alpha^{-1} \beta \sigma_j^2 \beta^{-1}\alpha) = 
(\alpha^{-1} \beta \sigma_j^2 \beta^{-1}\alpha) \sigma_i^2\\
\Leftrightarrow\quad&
\delta_i (\alpha^{-1} \beta \sigma_j^2 \beta^{-1}\alpha) = 
(\alpha^{-1} \beta \sigma_j^2 \beta^{-1}\alpha) \delta_i\\
\Leftrightarrow\quad&
(\beta^{-1} \alpha \delta_i \alpha^{-1} \beta) \sigma_j^2  = 
\sigma_j^2 (\beta^{-1} \alpha \delta_i \alpha^{-1} \beta)\\ 
\Leftrightarrow\quad&
(\beta^{-1} \alpha \delta_i \alpha^{-1} \beta) \delta_j = 
\delta_j (\beta^{-1} \alpha \delta_i \alpha^{-1} \beta)\\ 
\Leftrightarrow\quad&
(\alpha \delta_i \alpha^{-1}) (\beta \delta_j \beta^{-1}) =
(\beta \delta_j \beta^{-1}) (\alpha \delta_i \alpha^{-1})
\end{array}
\]
This shows that the bijection $\varphi\co  \hat \Upsilon \to \Upsilon$ extends to an isomorphism 
$\varphi\co  \MM(\hat \Omega)$ $\to \MM(\Omega)$. 
\end{proof}

Now, we have the following decomposition for $SB_n$.

\begin{lem}
$SB_n= \MM(\hat \Omega) \rtimes B_n$.
\end{lem}

\begin{proof}
Clearly, there exists a homomorphism $f\co  \MM(\hat \Omega) \rtimes B_n \to SB_n$ 
which sends $\beta$ to $\beta \in SB_n$ for all $\beta \in B_n$, and sends 
$\hat u$ to $\hat u \in SB_n$ for all $\hat u \in \hat \Upsilon$. On the other 
hand, one can easily verify using the presentation of $SB_n$ that there exists a homomorphism $g\co  
SB_n \to \MM(\hat \Omega) \rtimes B_n$ such that $g(\sigma_i^{\pm 1}) = \sigma_i^{\pm 1} \in 
B_n$ for all $i\in \{1, \dots, n-1\}$, and $g(\delta_i)=\delta_i \in \hat \Upsilon$ for all $i 
\in \{1, \dots, n-1\}$. Obviously, $f \circ g = \Id$ and $g \circ f = \Id$. 
\end{proof}

\begin{rmk}
\begin{enumerate}
\item Let $G(\hat \Omega)$ be the group given by the presentation
\[
G(\hat \Omega)= \langle \hat \Upsilon \ |\ \hat u \hat v = \hat v \hat u \qua \text{if}\qua \{\hat u, 
\hat v\} \in E(\hat \Omega) \rangle\,.
\]
It is well-known that $\MM(\hat\Omega)$ embeds in $G(\hat\Omega)$ (see \cite{DK}, \cite{Dub}), thus 
$SB_n= \MM(\hat \Omega) \rtimes B_n$ embeds in $G(\hat \Omega) \rtimes B_n$. This furnishes 
one more proof of the fact that $SB_n$ embeds in a group (see \cite{FKR}, \cite{Bas}, \cite{Key}).
\item The decomposition $SB_n= \MM(\hat \Omega) \rtimes B_n$ together with Lemma 2.2 can be used 
to solve the word problem in $SB_n$. The proof of this fact is left to the reader. Another 
solution to the word problem for $SB_n$ can be found in \cite{Cor}.
\end{enumerate}
\end{rmk}

\begin{proof}[Proof of Theorem 1.1] 
Consider the homomorphism $\deg\co  B_n \to \Z$ defined by 
$\deg(\sigma_i)=1$ for $1 \le i\le n-1$. For $k \in \Z$, let $B_n^{(k)}=\{ \beta \in 
B_n;\,\, \deg(\beta) =k\}$. We have the decomposition
\[
\Z[B_n]= \bigoplus_{k \in \Z} \Z[B_n^{(k)}]\,,
\]
where $\Z[B_n^{(k)}]$ denotes the free abelian group freely generated by $B_n^{(k)}$. Let $P 
\in \Z[B_n]$. We write $P=\sum_{k \in \Z} P_k$, where $P_k \in \Z[B_n^{(k)}]$ for all 
$k \in \Z$. Then $P_k$ is called the {\it $k$-th component} of $P$.

Let $\gamma, \gamma' \in SB_n$ such that $\eta(\gamma)= \eta(\gamma')$. We write $\gamma= \alpha 
\beta$ and $\gamma'= \alpha' \beta'$ where $\alpha, \alpha' \in \MM(\hat \Omega)$ and $\beta, 
\beta' \in B_n$ (see Lemma 2.3). Let $d=\deg(\beta)$. We observe that the $d$-th component 
of $\eta(\gamma)$ is $\pm\beta$, and, for $k<d$, the $k$-th component of $\eta(\gamma)$ is $0$. In particular, 
$\eta(\gamma)$ completely determines $\beta$. Since $\eta(\gamma)=\eta(\gamma')$, it 
follows that $\beta=\beta'$.

So, multiplying  $\gamma$ and $\gamma'$ on the right by $\beta^{-1}$ if necessary, we may assume 
that $\gamma=\alpha \in \MM(\hat \Omega)$ and $\gamma' = \alpha' \in \MM(\hat \Omega)$. 
Observe that
\[
(\nu \circ \varphi)(\gamma) = \eta(\gamma)= \eta(\gamma') = (\nu \circ \varphi) (\gamma')\,.
\]
Since $\nu$ is injective (Theorem 1.3) and $\varphi$ is an isomorphism (Lemma 2.2), we conclude 
that $\gamma=\gamma'$.
\end{proof}

\section{Proof of Theorem 1.3}

We start this section with the following result on graph monoids.

\begin{lem}
Let $\Gamma$ be a graph, let $X$ be the set of vertices, and let 
$E=E(\Gamma)$ be the set of edges of $\Gamma$. Let $x_1, \dots, x_l, y_1, \dots, y_l \in X$ and 
$k \in \{1,2, \dots, l\}$ such that:
\begin{itemize}
\item $x_1x_2 \dots x_l = y_1y_2 \dots y_l$ (in $\MM(\Gamma)$);
\item $y_k=x_1$, and $y_i\neq x_1$ for all $i=1, \dots, k-1$.
\end{itemize}
Then $\{y_i,x_1\}\in E(\Gamma)$ for all $i=1,2, \dots, k-1$.
\end{lem}

\begin{proof}
Let $F^+(X)$ denote the free monoid freely generated by $X$. Let $\equiv_1$ be the 
relation on $F^+(X)$ defined as follows. We set $u \equiv_1 v$ if there exist $u_1,u_2 \in 
F^+(X)$ and $x,y \in X$ such that $u=u_1xyu_2$, $v=u_1yxu_2$, and $\{x,y\} \in E(\Gamma)$. For $p 
\in \N$, we define the relation $\equiv_p$ on $F^+(X)$ by setting $u \equiv_p v$ if there 
exists a sequence $u_0=u, u_1, \dots, u_p=v$ in $F^+(X)$ such that $u_{i-1} \equiv_1 u_i$ for all 
$i=1, \dots, p$. Consider the elements $u=x_1x_2 \dots x_l$ and $v=y_1y_2 \dots y_l$ in $F^+(X)$. 
Obviously, there is some $p \in \N$ such that $u\equiv_p v$. Now, we prove the result of 
Lemma 3.1 by induction on $p$.

The case $p=0$ being obvious, we may assume $p\ge 1$. There exists a sequence $u_0=u, u_1, \dots, 
u_{p-1}, u_p=v$ in $F^+(X)$ such that $u_{i-1} \equiv_1 u_i$ for all $i=1, \dots, p$. By 
definition of $\equiv_1$, there exists $j \in \{1,2,\dots, l-1\}$ such that $\{y_j, y_{j+1}\} \in 
E(\Gamma)$ and $u_{p-1}= y_1 \dots y_{j-1} y_{j+1} y_j y_{j+2} \dots y_l$. If either $j<k-1$ or 
$j >k$, then, by the inductive hypothesis, we have $\{x_1,y_i\} \in E(\Gamma)$ for all $i=1, 
\dots, k-1$. If $j=k-1$, then, by the inductive hypothesis, we have $\{x_1,y_i\} \in E(\Gamma)$ 
for all $i=1, \dots, k-2$. Moreover, in this case, $\{y_j,y_{j+1}\}=\{y_{k-1},y_k\} = \{y_{k-
1},x_1\} \in E(\Gamma)$. If $j=k$, then, by the inductive hypothesis, we have $\{y_i,x_1\} \in 
E(\Gamma)$ for all $i=1, \dots, k-1$ and $i=k+1$.
\end{proof}

Now, consider the standard epimorphism $\theta\co  B_n \to \Sym_n$ defined by $\theta(\sigma_i) = 
(i,i+1)$ for $1 \le i\le n-1$. The kernel of $\theta$ is called the {\it pure braid group on $n$ 
strings}, and is denoted by $PB_n$. It has a presentation with generators
\[
A_{i\,j}= \sigma_{j-1} \dots \sigma_{i+1} \sigma_i^2 \sigma_{i+1}^{-1} \dots \sigma_{j-1}^{-1}\,, 
\quad 1 \le i<j\le n\,,
\]
and relations
\begin{gather*}
A_{r\,s}^{-1} A_{i\,j} A_{r\,s} = A_{i\,j}\quad 
\text{if}\qua r<s<i<j \text{ or } i<r<s<j\,,\\
A_{r\,s}^{-1} A_{i\,j} A_{r\,s} = A_{r\,j} A_{i\,j} A_{r\,j}^{-1}\quad 
\text{if}\qua s=i\,,\\
A_{r\,s}^{-1} A_{i\,j} A_{r\,s} = A_{i\,j} A_{s\,j} A_{i\,j} A_{s\,j}^{-1} A_{i\,j}^{-1}\quad 
\text{if}\qua i=r<s<j\,,\\
A_{r\,s}^{-1} A_{i\,j} A_{r\,s} = A_{r\,j} A_{s\,j} A_{r\,j}^{-1} A_{s\,j}^{-1} A_{i\,j} 
A_{s\,j} A_{r\,j} A_{s\,j}^{-1} A_{r\,j}^{-1} \quad
\text{if}\qua r<i<s<j\,.\\
\end{gather*}
(See \cite{Bir1}). We denote by $H_1(PB_n)$ the abelianization of $PB_n$, and, for $\beta \in PB_n$, 
we denote by $[\beta]$ the element of $H_1(PB_n)$ represented by $\beta$. A consequence of the 
above presentation is that $H_1(PB_n)$ is a free abelian group freely generated by $\{[A_{i\,j}];\,\, 
1 \le i<j\le n\}$. This last fact shall be of importance in the remainder of the paper.

For $1 \le i<j\le n$, we set
\[
\Upsilon_{i\,j}= \{\beta A_{i\,j} \beta^{-1}\ ;\ \beta \in PB_n\}\,.
\]

\begin{lem}
We have the disjoint union $\Upsilon=\bigsqcup_{i<j} \Upsilon_{i\,j}$.
\end{lem}

\begin{proof}
It is esily checked that
\[
\sigma_r A_{i\,j} \sigma_r^{-1}= \left\{
\begin{array}{ll}
A_{i\,j+1}\quad &\text{if}\qua r=j\,,\\
A_{j-1\,j} A_{i\,j-1} A_{j-1\,j}^{-1}\quad &\text{if}\qua r=j-1>i\,,\\
A_{i+1\,j}\quad &\text{if}\qua j-1>i=r\,,\\
A_{i\,j}^{-1} A_{i-1\,j} A_{i\,j}\quad&\text{if}\qua r=i-1\,,\\
A_{i\,j}\quad&\text{otherwise}\,.
\end{array}\right.
\]
This implies that the union
$\bigcup_{i<j} \Upsilon_{i\,j}$ is invariant by the action of $B_n$ by conjugation. Moreover, 
$\sigma_i^2=A_{i\,i+1} \in \Upsilon_{i\,i+1}$ for all $i \in \{1, \dots, n-1\}$, thus $\Upsilon \subset 
\bigcup_{i<j} \Upsilon_{i\,j}$. On the other hand, $A_{i\,j}$ is conjugate
(by an element of $B_n$) to $\sigma_i^2$, thus 
$\Upsilon_{i\,j} \subset \Upsilon$ for all $i<j$, therefore $\bigcup_{i<j} \Upsilon_{i\,j} 
\subset \Upsilon$.

Let $i,j,r,s \in \{1, \dots, n\}$ such that $i<j$, $r<s$, and $\{i,j\} \neq \{r,s\}$. Let $u \in 
\Upsilon_{i\,j}$ and $v \in \Upsilon_{r\,s}$. Then $[u]=[A_{i\,j}] \neq [A_{r\,s}] = [v]$, 
therefore $u \neq v$. 
This shows that $\Upsilon_{i\,j} \cap \Upsilon_{r\,s} = \emptyset$.
\end{proof}

The following lemmas 3.3 and 3.5 will be proved in Sections 4 and 5, respectively. 

Let $F(X)$ 
be a free group freely generated by some set $X$. Let $Y=\{g x g^{-1};\,\, g \in F(X)\ \text{and}\ 
x\in X\}$, and let $F^+(Y)$ be the free monoid freely generated by $Y$. We prove in Section 4 
that the homomorphism $\nu\co  F^+(Y) \to \Z[F(X)]$, defined by $\nu(y)=y-1$ for all $y \in Y$, 
is injective (Proposition 4.1). The proof of this result is based on the construction of a sort 
of Magnus expansion. Proposition 4.1 together with the fact that $PB_n$ can be decomposed as 
$PB_n= F \rtimes PB_{n-1}$, where $F$ is a free group freely generated by $\{A_{i\,n};\,\, 1\le i\le 
n-1\}$, are the main ingredients of the proof of Lemma 3.3. 

Choose some $x_0 \in X$, consider the decomposition $F(X)= \langle x_0 \rangle \ast 
F(X \setminus \{x_0\})$, and let $\rho\co  F(X) \to F(X)$ be an automorphism which fixes $x_0$ and which 
leaves $F(X \setminus \{x_0\})$ invariant. Let $y_1, \dots, y_l \in \{gx_0g^{-1};\,\, g \in F(X)\}$.
We prove in Section 5 that, if $\rho(y_1 \dots y_l)= y_1 \dots y_l$, then $\rho(y_i)=y_i$ for all
$i=1, \dots, l$ (Proposition 5.1). 
The proof of Lemma 3.5 is based 
on this result together with Corollary 3.4 below.

\begin{lem}
Let $i,j,r,s \in \{1, \dots, n\}$ such that $i<j$, $r<s$, $\{i,j\} \neq 
\{r,s\}$, and $\{i,j\} \cap \{r,s\} \neq \emptyset$. Let $\MM[i,j,r,s]$ be the free monoid 
freely generated by $\Upsilon_{i\,j} \cup \Upsilon_{r\,s}$, and let $\bar \nu\co  \MM[i,j,r,s] 
\to \Z[B_n]$ be the homomorphism defined by $\bar \nu(u)=u-1$ for all $u \in \Upsilon_{i\,j} 
\cup \Upsilon_{r\,s}$. Then $\bar \nu$ is injective.
\end{lem}

\begin{cor}
Let $i,j \in \{1, \dots, n\}$ such that $i<j$. Let $\MM[i,j]$ be the 
free monoid freely generated by $\Upsilon_{i\,j}$, and let $\bar \nu\co  \MM[i,j] \to 
\Z[B_n]$ be the homomorphism defined by $\bar\nu (u)=u-1$ for all $u \in \Upsilon_{i\,j}$. Then 
$\bar \nu$ is injective.
\end{cor}

\begin{lem}
Let $i,j,r,s\in \{1, \dots, n\}$ such that $i<j$, $r<s$, and $\{i,j\} \cap 
\{r,s\} = \emptyset$. (In particular, we have $n\ge 4$.) Let $\bar \Omega[i,j,r,s]$ be the graph 
defined as follows.
\begin{itemize}
\item $\Upsilon_{i\,j} \cup \Upsilon_{r\,s}$ is the set of vertices of $\bar\Omega 
[i,j,r,s]$;
\item $\{u,v\}$ is an edge of $\bar \Omega[i,j,r,s]$ if and only if we have $uv=vu$ in $B_n$.
\end{itemize}
Let $\MM[i,j,r,s]= \MM(\bar \Omega [i,j,r,s])$, and let $\bar \nu\co  \MM[i,j,r,s] \to  
\Z[B_n]$ be the homomorphism defined by $\bar \nu(u)=u-1$ for all $u \in \Upsilon_{i\,j} \cup 
\Upsilon_{r\,s}$. Then $\bar\nu$ is injective.
\end{lem}

\begin{proof}[Proof of Theorem 1.3] 
Recall the decomposition
\begin{equation}\label{eq31}
\Z[B_n]= \bigoplus_{k \in \Z} \Z[B_n^{(k)}]
\end{equation}
given in the proof of Theorem 1.1, where $B_n^{(k)}= \{\beta \in B_n;\,\, \deg(\beta)=k\}$, 
and $\Z [B_n^{(k)}]$ is the free abelian group freely generated by $B_n^{(k)}$. Note that 
$\deg(u)=2$ for all $u \in \Upsilon$.

Let $\alpha \in \MM(\Omega)$. We write $\alpha= u_1u_2 \dots u_l$, where $u_i \in \Upsilon$ 
for all $i=1, \dots, l$. Define the {\it length} of $\alpha$ to be $|\alpha|=l$. We denote by 
$\bar \alpha$ the element of $B_n$ represented by $\alpha$ (ie, $\bar\alpha= u_1 u_2 \dots u_l$ 
in $B_n$). Let $[1,l]=\{1,2,\dots, l\}$. Define a {\it subindex} of $[1,l]$ to be a sequence 
$I=(i_1,i_2, \dots, i_q)$ such that $i_1, i_2, \dots, i_q \in [1,l]$, and $i_1<i_2<\dots<i_q$. 
The notation $I\prec [1,l]$ means that $I$ is a subindex of $[1,l]$. The {\it length} of $I$ is 
$|I|=q$. For $I=(i_1,i_2, \dots, i_q) \prec [1,l]$, we set $\alpha(I)= u_{i_1} u_{i_2} \dots 
u_{i_q} \in \MM(\Omega)$
and $\bar\alpha(I)$ denotes the corresponding element of $B_n^{(2q)}$.

Observe that the decomposition of $\nu(\alpha)$ with respect to the direct sum \eqref{eq31} is:
\begin{equation}\label{eq32}
\nu(\alpha)= \sum_{q=0}^l (-1)^{l-q} \sum_{I \prec [1,l],\ |I|=q} \bar \alpha(I)\,,
\end{equation}
and
\[
\sum_{I \prec [1,l],\ |I|=q} \bar \alpha(I) \in \Z[B_n^{(2q)}]\,,
\]
for all $q=0,1, \dots, l$.

Let $\alpha' = u_1' u_2' \dots u_k' \in \MM(\Omega)$ such that $\nu(\alpha) = \nu(\alpha')$. 
The decomposition given in \eqref{eq32} shows that $k=l$ and
\begin{equation}\label{eq33}
\sum_{I \prec [1,l],\ |I|=q} \bar \alpha(I) = \sum_{I \prec [1,l],\ |I|=q} \bar \alpha'(I)\,, 
\end{equation}
for all $q=0,1, \dots, l$.

We prove that $\alpha = \alpha'$ by induction on $l$. The cases $l=0$ and $l=1$ being obvious, we 
assume $l \ge 2$.

Suppose first that $u_1'=u_1$. We prove
\begin{equation}\label{eq34}
\sum_{I \prec [2,l],\ |I|=q} \bar \alpha(I) = \sum_{I \prec [2,l],\ |I|=q} \bar \alpha'(I)
\end{equation}
by induction on $q$. The case $q=0$ being obvious, we assume $q \ge 1$. Then
\begin{align*}
&\sum_{I \prec [2,l],\ |I|=q} \bar \alpha(I)\\
=\ &\sum_{I \prec [1,l],\ |I|=q} \bar \alpha(I) -u_1 \cdot \sum_{I \prec [2,l],\ |I|=q-1} \bar 
\alpha(I)\\
=\ &\sum_{I \prec [1,l],\ |I|=q} \bar \alpha'(I) -u_1 \cdot \sum_{I \prec [2,l],\ |I|=q-1} \bar 
\alpha'(I) \quad \text{(by\ induction\ and\ \eqref{eq33})}\\
=\ &\sum_{I \prec [2,l],\ |I|=q} \bar \alpha'(I)\,.\\
\end{align*}
Let $\alpha_1= u_2 \dots u_l$ and $\alpha_1' = u_2' \dots u_l'$. By \eqref{eq34}, we have
\begin{align*}
\nu(\alpha_1) =&
\sum_{q=0}^{l-1} (-1)^{l-1-q} \sum_{I \prec [2,l],\ |I|=q} \bar \alpha(I)\\
=&\sum_{q=0}^{l-1} (-1)^{l-1-q} \sum_{I \prec [2,l],\ |I|=q} \bar \alpha'(I)=\nu(\alpha_1')
\end{align*}
thus, by the inductive hypothesis, $\alpha_1=\alpha_1'$, therefore $\alpha= u_1 \alpha_1 
= u_1 \alpha_1' = \alpha'$.

Now, we consider the general case. \eqref{eq33} applied to $q=1$ gives
\begin{equation}\label{eq35}
\sum_{i=1}^l u_i = \sum_{i=1}^l u_i' \,.
\end{equation}
So, there exists $k \in \{1, \dots, l\}$ such that $u_k'=u_1$ and $u_i' \neq u_1$ for all $i=1, 
\dots, k-1$. We prove that, for $1 \le i\le k-1$, $u_i'$ and $u_1=u_k'$ 
multiplicatively commute (in $B_n$ 
or, equivalently, in $\MM(\Omega)$). It follows that $\alpha'= u_1 u_1' \dots u_{k-1}' 
u_{k+1}' \dots u_l'$, and hence, by the case $u_1=u_1'$ considered before, $\alpha=\alpha'$.

Fix some $t \in \{1, \dots, k-1\}$. Let $i,j,r,s \in \{1, \dots, n\}$ such that $i<j$, $r<s$, 
$u_1=u_k' \in \Upsilon_{i\,j}$, and $u_t' \in \Upsilon_{r\,s}$. There are three possible cases 
that we handle simultaneously:
\begin{enumerate}
\item $\{i,j\} = \{r,s\}$;
\item $\{i,j\} \neq \{r,s\}$ and $\{i,j\} \cap \{r,s\} \neq \emptyset$;
\item $\{i,j\} \cap \{r,s\} = \emptyset$.
\end{enumerate}
Let $\bar \Omega [i,j,r,s]$ be the graph defined as follows. 
\begin{itemize}
\item $\Upsilon_{i\,j} \cup \Upsilon_{r\,s}$ is the set of vertices of $\bar \Omega 
[i,j,r,s]$;
\item $\{u,v\}$ is an edge of $\bar \Omega [i,j,r,s]$ if and only if we have $uv=vu$ in 
$B_n$.
\end{itemize}
Let $\MM[i,j,r,s]= \MM(\bar \Omega[i,j,r,s])$, and let $\bar \nu\co  \MM[i,j,r,s] \to 
\Z[B_n]$ be the homomorphism defined by $\bar \nu (u)=u-1$ for all $u \in \Upsilon_{i\,j} 
\cup \Upsilon_{r\,s}$. 
Note that, by Corollary 3.4 and Lemma 3.3, $\bar \Omega [i,j,r,s]$ has no edge and $\MM[i,j,r,s]$
is a free monoid in Cases 1 and 2. Moreover,
the homomorphism $\bar\nu$ is injective by Lemmas 3.3 and 3.5 and by 
Corollary 3.4.

Let $a_1=1,a_2, \dots, a_p \in [1,l]$, $a_1<a_2< \dots <a_p$, be the indices such that $u_{a_\xi} 
\in \Upsilon_{i\,j} \cup \Upsilon_{r\,s}$ for all $\xi=1,2, \dots, p$. Let $I_0=(a_1,a_2, \dots, 
a_p)$, and let $\alpha(I_0)= u_{a_1} u_{a_2} \dots u_{a_p} \in \MM[i,j,r,s]$. (It is true that 
$\MM[i,j,r,s]$ is a submonoid of $\MM(\Omega)$, but this fact is not needed for our 
purpose. So, we should consider $\alpha(I_0)$ as an element of $\MM[i,j,r,s]$, and not as an 
element of $\MM(\Omega)$.) Recall that, for $\beta \in PB_n$, we denote by $[\beta]$ the 
element of $H_1(PB_n)$ represented by $\beta$. Recall also that $H_1(PB_n)$ is a free abelian 
group freely generated by $\{[A_{i\,j}];\,\, 1\le i<j\le n\}$. Observe that
\begin{equation}\label{eq36}
\bar \nu(\alpha(I_0))= \sum_{q=0}^p (-1)^{p-q} \sum_{\substack{I \prec [1,l],\ |I|=q,\\ [\bar \alpha(I)] 
\in \Z[A_{i\,j}] + \Z[A_{r\,s}]}} \bar \alpha(I)\,.
\end{equation}
Let $b_1, \dots, b_p \in [1,l]$, $b_1<b_2< \dots<b_p$, be the indices such that $u_{b_\xi}' \in 
\Upsilon_{i\,j} \cup \Upsilon_{r\,s}$ for all $\xi=1,2, \dots, p$. (Clearly, \eqref{eq35} implies that we have 
as many $a_\xi$'s as $b_\xi$'s.) Note that $t,k \in \{b_1, \dots, b_p\}$. Let $I_0'=(b_1, b_2, 
\dots, b_p)$, and let $\alpha'(I_0')=u_{b_1}' u_{b_2}' \dots u_{b_p}' \in \MM[i,j,r,s]$. By \eqref{eq33}
we have
\[
\sum_{\substack{ I \prec [1,l],\ |I|=q,\\ [\bar \alpha(I)] \in \Z[A_{i\,j}] + \Z[A_{r\,s}]}} 
\bar \alpha(I) = \sum_{\substack{I \prec [1,l],\ |I|=q,\\ [\bar \alpha'(I)] \in \Z[A_{i\,j}] +  
\Z[A_{r\,s}]}} \bar \alpha'(I)\,,
\]
for all $q \in \N$, thus, by \eqref{eq36}, $\bar \nu (\alpha(I_0))= \bar \nu (\alpha'(I_0'))$. 
Since $\bar\nu$ is injective, it follows that $\alpha(I_0)=\alpha'(I_0')$, and we conclude by Lemma 
3.1 that $u_t'$ and $u_k'=u_1$ commute.
\end{proof}

\section{Proof of Lemma 3.3}

As pointed out in the previous section, the key point of the proof of Lemma 3.3 is the following result.

\begin{prop}
Let $F(X)$ be a free group freely generated by some set $X$, let 
$Y=\{ gxg^{-1};\,\, g \in F(X)\ \text{and}\ x \in X\}$, let $F^+(Y)$ be the free monoid freely 
generated by $Y$, and let $\nu\co  F^+(Y) \to \Bbb Z[F(X)]$ be the homomorphism defined by 
$\nu(y)=y-1$ for all $y \in Y$. Then $\nu$ is injective.
\end{prop}

First, we shall prove Lemmas 4.2, 4.3, and 4.4 that are preliminary results to the proof of 
Proposition 4.1.

Let $\deg\co  F(X) \to \Z$ be the homomorphism defined by $\deg(x)=1$ for all $x 
\in X$. Write $\AA=\Z[F(X)]$. For $k \in \Z$, let $F_k(X)=\{ g\in F(X);\,\, \deg(g) 
\ge k\}$, and let $\AA_k=\Z[F_k(X)]$ be the free $\Z$--module freely generated by 
$F_k(X)$. The family $\{ \AA_k\}_{k \in \Z}$ is a filtration of $\AA$ compatible with 
the multiplication, that is:
\begin{itemize}
\item $\AA_k \subset \AA_l$ if $k \ge l$;
\item $\AA_p \cdot \AA_q \subset \AA_{p+q}$ for all $p,q \in \Z$;
\item $1 \in \AA_0$.
\end{itemize}
Moreover, this filtration is a separating filtration, that is:
\begin{itemize}
\item $\cap_{k \in \Z} \AA_k = \{0\}$.
\end{itemize}
Let $\tilde\AA$ denote the completion of $\AA$ with respect to this filtration. For $k \in 
\Z$, we write $F^{(k)}(X)= \{g\in F(X);\,\, \deg(g)=k\}$, and we denote by $\AA^{(k)}= 
\Z[F^{(k)}(X)]$ the free $\Z$--module freely generated by $F^{(k)}(X)$. Then any element 
of $\tilde \AA$ can be uniquely represented by a formal series $\sum_{k=d}^{+\infty} P_k$, 
where $d \in \Z$ and $P_k \in \AA^{(k)}$ for all $k \ge d$.

We take a copy $G_x$ of $\Z \times \Z$ generated by $\{x,\hat x\}$, for all $x \in X$, and 
we set $\hat G=\ast_{x\in X} G_x$. Let $\UU(\tilde \AA)$ denote the group of units of 
$\tilde \AA$. Then there is a homomorphism $\hat \eta\co  \hat G \to \UU(\tilde\AA)$ 
defined by
\[
\hat \eta(x)=x\,,\quad \hat \eta(\hat x)= x-1\,, \quad\text{for}\ x\in X\,.
\]
Note that
\[
\hat \eta(\hat x^{-1})= -\sum_{k=0}^{+\infty} x^k\,, \quad\text{for}\ x\in X\,.
\]
The homomorphism $\hat \eta$ defined above is a sort of Magnus expansion and the proof of the 
following lemma is strongly inspired by the proof of \cite{Bou}, Ch. II, \S\ 5, Thm. 1.

\begin{lem}
The homomorphism $\hat \eta\co  \hat G \to \UU(\tilde \AA)$ is 
injective.
\end{lem}

\begin{proof}
Let $g \in \hat G$. Define the {\it normal form} of $g$ to be the finite sequence 
$(g_1,g_2, \dots, g_l)$ such that:
\begin{itemize}
\item for all $i \in \{1, \dots, l\}$, there exists $x_i \in X$ such that $g_i \in G_{x_i} 
\setminus \{1\}$;
\item $x_i \neq x_{i+1}$ for all $i=1, \dots, l-1$;
\item $g=g_1g_2 \dots g_l$.
\end{itemize}
Clearly, such an expression for $g$ always exists and is unique. The {\it length} of $g$ is 
defined to be $\lg(g)=l$.

Let $(p,q)\in \Z \times \Z$, $(p,q) \neq (0,0)$. Write
\[
(t-1)^pt^q= \sum_{k=d}^{+\infty} c_{k\,p\,q} t^k\,,
\]
where $d \in \Z$ and $c_{k\,p\,q} \in \Z$ for all $k \ge d$. We show that there exists 
$a \ge d$ such that $a \neq 0$ and $c_{a\,p\,q} \neq 0$. If $q \neq 0$, then $a=q\neq 0$ and 
$c_{q\,p\,q}= \pm 1 \neq 0$. If $q=0$, then $a=1 \neq 0$ and $c_{1\,p\,0}=\pm p\neq 0$.

Let $g \in \hat G$, $g \neq 1$. Let $(\hat x_1^{p_1} x_1^{q_1}, \dots, \hat x_l^{p_l} x_l^{q_l})$ 
be the normal form of $g$. We have
\begin{align*}
\hat\eta (g)&= (x_1-1)^{p_1} x_1^{q_1} (x_2-1)^{p_2} x_2^{q_2} \dots (x_l-1)^{p_l} x_l^{q_l}\\
&= \sum_{k_1\ge d_1, \dots,k_l\ge d_l} 
c_{k_1\,p_1\,q_1} c_{k_2\,p_2\,q_2} \dots c_{k_l\,p_l\,q_l} \cdot x_1^{k_1} 
x_2^{k_2} \dots x_l^{k_l}\,.
\end{align*}
By the above observation, there exist $a_1,a_2, \dots, a_l \in \Z \setminus \{0\}$ such that 
$c_{a_i\,p_i\,q_i}\neq 0$ for all $i=1, \dots, l$. Now, we show that $x_1^{k_1} \dots x_l^{k_l} 
\neq x_1^{a_1} \dots x_l^{a_l}$ if $(k_1, \dots, k_l)\neq (a_1, \dots, a_l)$. This implies that 
the coefficient of $x_1^{a_1} \dots x_l^{a_l}$ in $\hat\eta(g)$ is \break
$c_{a_1\,p_1\,q_1} \dots 
c_{a_l\,p_l\,q_l}\neq 0$, thus $\hat\eta(g) \neq 1$.

Since $(\hat x_1^{p_1} x_1^{q_1}, \dots, \hat x_l^{p_l} x_l^{q_l})$ is the normal form of $g$, we 
have $x_i \neq x_{i+1}$ for all $i=1, \dots, l-1$, thus $(x_1^{a_1}, \dots, x_l^{a_l})$ is the 
normal form of $x_1^{a_1} \dots x_l^{a_l}$. Suppose $k_i \neq 0$ for all $i=1, \dots, l$. Then 
$(x_1^{k_1}, \dots, x_l^{k_l})$ is the normal form of $x_1^{k_1} \dots x_l^{k_l}$, therefore 
$x_1^{k_1} \dots x_l^{k_l} \neq x_1^{a_1} \dots x_l^{a_l}$ if $(k_1, \dots, k_l) \neq (a_1, 
\dots, a_l)$. Suppose there exists $i \in \{1, \dots, l\}$ such that $k_i=0$. Then
\[
\lg (x_1^{k_1} \dots x_l^{k_l}) <l=\lg(x_1^{a_1} \dots x_l^{a_l})\,,
\]
thus $x_1^{k_1} \dots x_l^{k_l} \neq x_1^{a_1} \dots x_l^{a_l}$.
\end{proof}

For each $x \in X$, we take a copy $SG_x$ of $\Z \times \N$ generated as a monoid by 
$\{x,x^{-1}, \hat x\}$, and we set $SG=\ast_{x \in X} SG_x$. Then there is a homomorphism $\eta\co  
SG \to \Z[F(X)]$ defined by
\[
\eta(x^{\pm 1}) = x^{\pm 1}\,, \quad \eta(\hat x)=x-1\,, \quad \text{for}\ x\in X\,.
\]

\begin{lem}
The homomorphism $\eta\co  SG \to \Z[F(X)]$ is injective.
\end{lem}

\begin{proof}
We have $SG \subset \hat G$, and, since $\{ \AA_k\}_{k \in \Z}$ is a separating 
filtration, $\AA=\Z[F(X)]$ is a subalgebra of $\tilde \AA$. Now, observe that $\eta\co  SG 
\to \Z[F(X)]$ is the restriction of $\hat \eta$ to $SG$, thus, by Lemma 4.2, $\eta$ is 
injective. 
\end{proof}

Let $\hat Y=\{g \hat x g^{-1};\,\, g \in F(X)\ \text{and}\ x\in X\} \subset SG$, and let $F^+(\hat 
Y)$ be the free monoid freely generated by $\hat Y$. The proof of the following lemma is left to 
the reader.
A more general statement can be found in \cite{DK}.

\begin{lem}
We have $SG=F^+(\hat Y) \rtimes F(X)$.
\end{lem}

Now, we can prove Proposition 4.1, and, consequently, Lemma 3.3.

\begin{proof}[Proof of Proposition 4.1]
Let $\hat \nu\co  F^+(\hat Y) \to \Z[F(X)]$ be the restriction 
of $\eta\co  SG= F^+(\hat Y) \rtimes F(X) \to \Z[F(X)]$ to $F^+(\hat Y)$, and let $\varphi\co  
F^+(\hat Y) \to F^+(Y)$ be the epimorphism defined by $\varphi(g \hat x g^{-1}) = gxg^{-1}$ for 
all $g \in F(X)$ and $x \in X$. (The proof that $\varphi$ is well-defined is left to the reader.) 
The homomorphism $\hat \nu$ is injective (Lemma 4.3), $\varphi$ is a surjection, and $\hat \nu = 
\nu \circ \varphi$, thus $\varphi$ is an isomorphism and $\nu$ is injective. 
\end{proof}

\begin{proof}[Proof of Lemma 3.3] 
Take $\zeta \in \Sym_n$ such that $\zeta(\{i,j\}) = \{1,n\}$ and \break 
$\zeta(\{r,s\}) = \{2,n\}$. Choose $\beta \in B_n$ such that $\theta (\beta)= \zeta$. Then 
$\beta \Upsilon_{i\,j} \beta^{-1} = \Upsilon_{1\,n}$ and $\beta \Upsilon_{r\,s} \beta^{-1} = 
\Upsilon_{2\,n}$. So, up to conjugation by $\beta$ if necessary, we may assume that 
$\{i,j\}=\{1,n\}$ and $\{r,s\}=\{2,n\}$.

Let $F$ be the subgroup of $PB_n$ generated by $\{A_{i\,n};\,\, 1\le i\le n-1\}$. We have:
\begin{enumerate}
\item $F$ is a free group freely generated by $\{A_{i\,n};\,\, 1\le i\le n-1\}$;
\item $PB_n=F \rtimes PB_{n-1}$;
\item $\Upsilon_{i\,n}= \{g A_{i\,n} g^{-1};\,\, g\in F\}$ for all $i=1, \dots, n-1$.
\end{enumerate}
(1) and (2) are well-known and are direct consequences of the presentation of $PB_n$ given in 
Section 3, and  (3) follows from the fact that the conjugacy class of $A_{i\,n}$ in $F$ is 
invariant by the action of $PB_{n-1}$.

Let $\Upsilon'= \sqcup_{i=1}^{n-1} \Upsilon_{i\,n}$, and let $F^+(\Upsilon')$ be the free monoid 
freely generated by $\Upsilon'$. By Proposition 4.1, the homomorphism $\nu'\co  F^+(\Upsilon') \to 
\Z[F]$, defined by $\nu'(u)=u-1$ for all $u \in \Upsilon'$, is injective. Recall that  
$\MM[1,n,2,n]$ denotes the free monoid freely generated by $\Upsilon_{1\,n} \cup \Upsilon_{2\,n}$. 
Then $\MM[1,n,2,n] \subset F^+(\Upsilon')$, $\Z[F] \subset \Z[B_n]$, and $\bar\nu\co  
\MM[1,n,2,n] \to \Z[B_n]$ is the restriction of $\nu'$ to $\MM[1,n,2,n]$, thus 
$\bar\nu$ is injective. 
\end{proof}

\section{Proof of Lemma 3.5}

We assume throughout this section that $n \ge 4$. As pointed out in Section 3, one of the main 
ingredients of the proof of Lemma 3.5 is the following result.

\begin{prop}
Let $F(X)$ be a free group freely generated by some set $X$, let $x_0 
\in X$, and let $\rho\co  F(X) \to F(X)$ be an automorphism 
which fixes $x_0$ and leaves $F(X \setminus \{x_0\})$ invariant
(where $F(X \setminus \{x_0\})$ denotes the subgroup 
of $F(X)$ (freely) generated by $X \setminus \{x_0\}$). Let $y_1, \dots, y_l \in \{ gx_0g^{-1};\,\, 
g\in F(X)\}$. If $\rho(y_1y_2 \dots y_l)= y_1y_2 \dots y_l$, then $\rho(y_i)=y_i$ for all $i=1, 
\dots, l$.
\end{prop}

\proof
Let $Z=\{hx_0h^{-1};\,\, h\in F(X \setminus \{x_0\}) \}$, and let $F(Z)$ be the subgroup 
of $F(X)$ generated by $Z$. Observe that $Z$ freely generates $F(Z)$, $\rho$ permutes the 
elements of $Z$, and $\{gx_0g^{-1};\,\, g\in F(X)\} = \{ \beta z \beta^{-1};\,\, \beta \in F(Z)\ 
\text{and}\ z\in Z\}$.

For $f \in F(Z)$, we denote by $\lg (f)$ the word length of $f$ with respect to $Z$. For 
$f,g \in F(Z)$, we write $fg = f \ast g$ if $\lg(fg)= \lg(f) + \lg(g)$. Note 
that, if $fg=f \ast g$, then $\rho(fg)= \rho(f) \ast \rho(g)$. Moreover, if $fg=f \ast g$ and 
$\rho(fg)=fg$, then $\rho(f)=f$ and $\rho(g)=g$.

Let $g_0=y_1y_2 \dots y_l$. Recall that we are under the assumption that $\rho(g_0)=g_0$. For 
$i=1, \dots, l$, let $\beta_i \in F(Z)$ and $z_i \in Z$ such that $y_i= \beta_i \ast z_i \ast 
\beta_i^{-1}$. Now, we prove that $\rho(y_i)=y_i$ for all $i=1, \dots, l$ by induction on 
$\sum_{i=1}^l \lg(y_i) = l + 2\sum_{i=1}^l \lg(\beta_i)$.

We have three cases to study.

\begin{case1}
There exists $t \in \{1, \dots, l-1\}$ such that $\beta_{t+1}= \beta_t \ast z_t^{-
1} \ast \gamma_t$, where $\gamma_t \in F(Z)$.

Let $y_{t+1}'= \beta_t\gamma_t z_{t+1} \gamma_t^{-1} \beta_t^{-1} = y_t y_{t+1} y_t^{-1}$. 
Observe that
\[
g_0=y_1 \dots y_{t-1} y_{t+1}' y_t y_{t+2} \dots y_l\,.
\]
We have $\lg( \beta_t \gamma_t) < \lg (\beta_{t+1})$, thus, by the inductive 
hypothesis, $\rho(y_i)= y_i$ for all $i=1, \dots, t-1,t, t+2, \dots, l$, and 
$\rho(y_{t+1}')=y_{t+1}'$. Moreover, since $y_{t+1}= y_t^{-1} y_{t+1}' y_t$, we also have 
$\rho(y_{t+1})=y_{t+1}$.
\end{case1}

\begin{case2}
There exists $t \in \{2, \dots, l\}$ such that $\beta_{t-1} = \beta_t \ast z_t \ast 
\gamma_t$, where $\gamma_t \in F(Z)$.

Then we prove that $\rho(y_i)=y_i$ for all $i=1, \dots, l$ as in the previous case.
\end{case2}

\begin{case3}
For all $t \in \{1, \dots, l\}$ and for all $\gamma_t \in F(Z)$ we have 
$\beta_{t+1} \neq \beta_t \ast z_t^{-1} \ast \gamma_t$ and $\beta_{t-1} \neq \beta_t \ast z_t 
\ast \gamma_t$.

We observe that
\[
g_0= \beta_1 \ast z_1 \ast \beta_1^{-1}\beta_2 \ast z_2 \ast \dots \ast \beta_{l-1}^{-1} \beta_l \ast z_l 
\ast \beta_l^{-1} \,.
\]
Since $\rho(g_0)=g_0$, it follows that $\rho(\beta_1)=\beta_1$, $\rho(z_1)=z_1$, $\rho(\beta_1^{-
1} \beta_2)= \beta_1^{-1} \beta_2$, $\rho(z_2)=z_2$, \dots, $\rho(\beta_{l-1}^{-1} \beta_l)= 
\beta_{l-1}^{-1} \beta_l$, $\rho(z_l)=z_l$, and $\rho(\beta_l^{-1})= \beta_l^{-1}$. This clearly 
implies that $\rho(y_i)=y_i$ for all $i=1, \dots, l$.\qed
\end{case3}

\begin{cor}
Let $u\in \Upsilon_{1\,2}$ and $v_1, \dots, v_l \in \Upsilon_{n-1\,n}$. 
If $u$ commutes with $v_1v_2 \dots v_l$ (in $B_n$), then $u$ commutes with $v_i$ for all $i=1, 
\dots, l$.
\end{cor}

\begin{proof}
Let $\alpha_0 \in PB_n$ such that $u= \alpha_0 A_{1\,2} \alpha_0^{-1}$. Up to 
conjugation of $v_1, \dots, v_l$ by $\alpha_0^{-1}$ if necessary, we can suppose that 
$\alpha_0=1$ and $u=A_{1\,2}$.

Recall that $F$ denotes the subgroup of $PB_n$ generated by $\{A_{i\,n};\,\, 1\le i\le n-1\}$. Recall 
also that:
\begin{itemize}
\item $F$ is a free group freely generated by $\{ A_{i\,n};\,\, 1\le i\le n-1\}$;
\item $PB_n= F \rtimes PB_{n-1}$;
\item $\Upsilon_{i\,n}= \{gA_{i\,n} g^{-1};\,\, g\in F\}$ for all $i=1, \dots, n-1$.
\end{itemize}
Let $\rho\co  F \to F$ be the action of $A_{1\,2}$ by conjugation on $F$ (namely, $\rho(g)= A_{1\,2} 
g A_{1\,2}^{-1}$). Observe that $\rho(A_{n-1\,n})= A_{n-1\,n}$ and the subgroup of $F$ generated 
by $\{A_{i\,n};\,\, 1\le i\le n-2\}$ is invariant by $\rho$. Then, Proposition 5.1 shows that 
$\rho(v_i)=v_i$ for all $i=1, \dots, l$ if $\rho(v_1v_2 \dots v_l)= v_1v_2 \dots v_l$. 
\end{proof}
 
Now, we can prove Lemma 3.5.

\begin{proof}[Proof of Lemma 3.5] 
Take $\zeta \in \Sym_n$ such that $\zeta(\{i,j\}) = \{1,2\}$ and \break
$\zeta( \{r,s\}) = \{n-1,n\}$. Choose $\beta \in B_n$ such that $\theta (\beta)= \zeta$. Then 
$\beta \Upsilon_{i\,j} \beta^{-1} = \Upsilon_{1\,2}$ and $\beta \Upsilon_{r\,s} \beta^{-1} = 
\Upsilon_{n-1\,n}$. So, up to conjugation by $\beta$ if necessary, we may assume that 
$\{i,j\}=\{1,2\}$ and $\{r,s\}=\{n-1,n\}$.

We use the same notations as in the proof of Theorem 1.3. Let $\alpha \in \MM[1,2,n-1,n]$. We 
write $\alpha= u_1u_2 \dots u_l$, where $u_i \in \Upsilon_{1\,2} \cup \Upsilon_{n-1\,n}$ for all 
$i=1, \dots, l$. Define the {\it length} of $\alpha$ to be $|\alpha|=l$. We denote by $\bar 
\alpha$ the element of $B_n$ represented by $\alpha$. Let $[1,l]=\{1, 2, \dots, l\}$. Define a 
{\it subindex} of $[1,l]$ to be a sequence $I=(i_1, i_2, \dots, i_q)$ such that $i_1,i_2, \dots, 
i_q \in [1,l]$ and $i_1<i_2< \dots <i_q$. The notation $I\prec [1,l]$ means that $I$ is a 
subindex of $[1,l]$. The {\it length} of $I$ is $|I|=q$. For $I= (i_1,i_2, \dots, i_q) \prec 
[1,l]$, we set $\alpha(I)= u_{i_1} u_{i_2} \dots u_{i_q} \in \MM[1,2,n-1,n]$.

Observe that
\begin{equation}\label{eq51}
\bar \nu (\alpha) = \sum_{q=0}^l (-1)^{l-q} \sum_{I\prec [1,l],\ |I|=q} \bar \alpha(I)\,, 
\end{equation}
and
\[
\sum_{I\prec [1,l],\ |I|=q} \bar \alpha(I) \in \Bbb Z [B_n^{(2q)}]\,,
\]
for all $q=0,1, \dots, l$.

Let $\alpha' = u_1' u_2' \dots u_k' \in \MM[1,2,n-1,n]$ such that $\bar \nu(\alpha) = \bar \nu 
(\alpha')$. As in the proof of Theorem 1.3, the decomposition given in \eqref{eq51} shows that $k=l$ and
\begin{equation}\label{eq52}
\sum_{I\prec [1,l],\ |I|=q} \bar \alpha(I) = \sum_{I\prec [1,l],\ |I|=q} \bar \alpha'(I)\,,
\end{equation}
for all $q=0,1, \dots, l$.

We prove that $\alpha=\alpha'$ by induction on $l$. The cases $l=0$ and $l=1$ being obvious, we 
assume $l \ge 2$. 

Assume first that $u_1'=u_1$. Then, by the same argument as in the proof of Theorem~1.3, 
$\alpha=\alpha'$.

Now, we consider the general case. \eqref{eq52} applied to $q=1$ gives
\begin{equation}\label{eq53}
\sum_{i=1}^l u_i = \sum_{i=1}^l u_i'\,.
\end{equation}
It follows that there exists a permutation $\zeta \in \Sym_l$ such that $u_i = u_{\zeta 
(i)}'$ for all $i=1, \dots,l$. (Note that the permutation $\zeta \in \Sym_l$ is not 
necessarily unique. Actually, $\zeta$ is unique if and only if $u_i \neq u_j$ for all $i \neq 
j$.)

Let $a_1,a_2, \dots, a_p \in [1,l]$, $a_1<a_2< \dots <a_p$, be the indices such that $u_{a_\xi} 
\in \Upsilon_{1\,2}$ for all $\xi=1, \dots, p$. Let $I_0=(a_1,a_2, \dots, a_p)$. Recall that, for 
$\beta \in PB_n$, we denote by $[\beta]$ the element of $H_1(PB_n)$ represented by $\beta$. 
Recall also that $H_1(PB_n)$ is a free abelian group freely generated by $\{[A_{i\,j}];\,\, 1 \le 
i<j\le n\}$. Observe that $\alpha(I_0) \in \MM[1,2]$ and
\begin{equation}\label{eq54}
\bar \nu (\alpha (I_0)) = \sum_{k=0}^p (-1)^{p-k} \sum_{\substack{I \prec [1,l],\ |I|=k,\\ [\bar 
\alpha(I)] \in \Z[A_{1\,2}]}} \bar \alpha(I)\,.
\end{equation}
Let $a_1',a_2', \dots, a_p' \in [1,l]$, $a_1'<a_2'< \dots <a_p'$, be the indices such that 
$u_{a_\xi'}' \in \Upsilon_{1\,2}$ for all $\xi=1, \dots, p$. Note that $\{ \zeta(a_1'), 
\zeta(a_2'), \dots, \zeta(a_p') \} = \{a_1,a_2, \dots, a_p\}$. Let $I_0'=(a_1',a_2', \dots, 
a_p')$. By \eqref{eq52}, we have
\begin{equation}\label{eq55}
\sum_{\substack{I \prec [1,l],\ |I|=k,\\ [\bar \alpha(I)] \in \Z[A_{1\,2}]}} \bar \alpha(I) = 
\sum_{\substack{I \prec [1,l],\ |I|=k,\\ [\bar \alpha'(I)] \in \Z[A_{1\,2}]}} \bar 
\alpha'(I)\,,
\end{equation}
for all $k \in \N$, thus, by \eqref{eq54}, $\bar \nu (\alpha(I_0)) = \bar \nu(\alpha'(I_0'))$. By 
Corollary 3.4, it follows that $\alpha(I_0)=\alpha'(I_0')$. So, $u_{a_i'}'=u_{a_i}$ for all 
$i=1, \dots, p$, and
the permutation $\zeta \in 
\Sym_l$ can be chosen so that $\zeta (a_i')=a_i$ for all $i=1, \dots, p$.

Let $b_1,b_2, \dots, b_q \in [1,l]$, $b_1<b_2< \dots <b_q$, be the indices such that $u_{b_\xi} 
\in \Upsilon_{n-1\,n}$ for all $\xi=1, \dots, q$. Note that $[1,l]= \{a_1, \dots, a_p, b_1, 
\dots, b_q\}$. Let $J_0=(b_1, b_2, \dots, b_q)$. Let $b_1',b_2', \dots, b_q' \in [1,l]$, $b_1'< 
b_2'< \dots < b_q'$, be the indices such that $u_{b_\xi'}' \in \Upsilon_{n-1\,n}$ for all 
$\xi=1, \dots, q$, and let $J_0'=(b_1',b_2', \dots, b_q')$. We also have $\alpha(J_0)= 
\alpha'(J_0') \in \MM[n-1,n]$, $u_{b_i}=u_{b_i'}'$ for all $i=1, \dots, q$, and $\zeta$ can be 
chosen so that $\zeta (b_i')=b_i$ for all $i=1, \dots, q$.

Without loss of generality, we can assume that $u_1 \in \Upsilon_{1\,2}$ (namely, $a_1=1$). Let 
$i \in \{1, \dots, p\}$. We set:
\begin{align*}
S(i)=&\left\{\begin{array}{ll}
0&\text{if}\qua a_i<b_1\,,\\
j&\text{if}\qua b_j<a_i<b_{j+1}\,,\\
q&\text{if}\qua b_q<a_i\,.
\end{array}\right.\\
T(i)=&\left\{\begin{array}{ll}
0&\text{if}\qua a_i'<b_1'\,,\\
j&\text{if}\qua b_j'<a_i'<b_{j+1}'\,,\\
q&\text{if}\qua b_q'<a_i'\,.
\end{array}\right.
\end{align*}
Note that $\alpha'=u_{b_1'}' \dots u_{b_{T(1)}'}' u_{a_1'}' \dots = u_{b_1} \dots u_{b_{T(1)}} 
u_{a_1} \dots$. Now, we show that $u_1=u_{a_1}$ commutes with $u_{b_i}$ for all $i=1, \dots, 
T(1)$. It follows that $\alpha'=u_1 u_{b_1} \dots u_{b_{T(1)}} \dots$, and hence, by the case 
$u_1'=u_1$ considered before, $\alpha=\alpha'$.

Let
\begin{align*}
v_i=& u_{b_1} \dots u_{b_{S(i)}} u_{a_i} u_{b_{S(i)}}^{-1} \dots u_{b_1}^{-1} \in 
\Upsilon_{1\,2}\,,\\
v_i'=& u_{b_1} \dots u_{b_{T(i)}} u_{a_i} u_{b_{T(i)}}^{-1} \dots 
u_{b_1}^{-1} \in \Upsilon_{1\,2}\,,
\end{align*}
for all $i=1, \dots, p$, and let
\[
\gamma= v_1v_2 \dots v_p \in \MM[1,2]\,, \quad \gamma'= v_1'v_2' \dots v_p' \in \MM[1,2]\,.
\]
Observe that
\begin{align*}
\bar \nu (\gamma) &= \left( \sum_{k=0}^p (-1)^{p-k} \sum_{\substack{I \prec [1,l],\ |I|=k+q,\\ [\bar 
\alpha(I)] = k[A_{1\,2}] +q[A_{n-1\,n}]}} \bar \alpha(I) \right) \bar \alpha(J_0)^{-1}\,,\\
\bar \nu (\gamma') &= \left( \sum_{k=0}^p (-1)^{p-k} \sum_{\substack{I \prec [1,l],\ |I|=k+q,\\ [\bar 
\alpha'(I)] = k[A_{1\,2}] +q[A_{n-1\,n}]}} \bar \alpha'(I) \right) \bar \alpha'(J_0')^{-
1}\,.
\end{align*}
We know that $\alpha(J_0) = \alpha'(J_0')$, and, by \eqref{eq52}, 
\[
\sum_{\substack{I \prec [1,l],\ |I|=k+q,\\ [\bar \alpha(I)] = k[A_{1\,2}] +q[A_{n-1\,n}]}} \bar 
\alpha(I) = \sum_{\substack{I \prec [1,l],\ |I|=k+q,\\ [\bar \alpha'(I)] = k[A_{1\,2}] +q[A_{n-1\,n}]}} 
\bar \alpha'(I)\,,
\]
for all $k=0,1, \dots, p$, thus $\bar \nu (\gamma) = \bar \nu (\gamma')$. By Corollary 3.4, it 
follows that $\gamma=\gamma'$, namely, $v_i=v_i'$ for all $i=1, \dots, p$. So,
\[
u_1=v_1=v_1'= u_{b_1} \dots u_{b_{T(1)}} u_1 u_{b_{T(1)}}^{-1} \dots u_{b_1}^{-1}\,,
\]
thus $u_1$ and $u_{b_1} \dots u_{b_{T(1)}}$ commute (in $B_n$). We conclude by Corollary 5.2 that 
$u_1$ and $u_{b_i}$ commute for all $i=1, \dots, T(1)$. 
\end{proof}

\bibliographystyle{gtart}

\end{document}